\documentclass{amsart}
 \usepackage{amsaddr}
\usepackage[utf8]{inputenc}
\usepackage[a4paper, bottom = 33mm, top = 33mm]{geometry}
\usepackage{url}
\usepackage{color}
\usepackage{amsmath, amsthm, amssymb, amsfonts, amscd}
\usepackage{graphicx}
\usepackage{graphicx,caption}
\usepackage{multicol}

\newcommand{\XX}{\mathbb{X}}
\newcommand{\OO}{\mathbb{O}}

\title{GridPyM: a Python module to handle grid diagrams}
\author{Agnese Barbensi and Daniele Celoria}
\address{School of Mathematics and Statistics, University of Melbourne}
\date{}

\begin{document}

\maketitle
\begin{abstract}
Grid diagrams are a combinatorial version of classical link diagrams, widely used in theoretical, computational and applied knot theory. Motivated by questions from (bio)-physical knot theory, we introduce GridPyM, a Sage compatible Python module that handles grid diagrams. GridPyM focuses on generating and simplifying grids, and on modelling local transformations between them.
\end{abstract}

\section{Introduction}

Grid diagrams are a handy and computationally efficient model for oriented link diagrams. 
They admit a concise description, can be easily randomised and --unlike standard link diagrams-- are able to encode certain local geometric features of links. They were first introduced in a slightly different form in~\cite{lyon}, and later on in~\cite{dynnikov,CROMWbasic} under the name of \emph{arc presentations}. The latter paper also introduced a set of combinatorial moves, commonly known as Cromwell moves, connecting any two grids representing equivalent link diagrams.

Grid diagrams are ubiquitous in theoretical and applied knot theory. Indeed, they arise naturally in contact geometry, as they represent Legendrian diagrams of links in $S^3$, endowed with its standard tight contact structure. Further, a certain subset of Cromwell moves bijectively corresponds to Reidemeister moves preserving the Legendrian isotopy class of the link~\cite{geiges,NG}.

Grid diagrams can also be interpreted as representing links in a genus one Heegaard decomposition of the three sphere. Here, the grid is identified with the Heegaard torus, with horizontal and vertical lines acting as the cores and co-cores of the $1$ and $2$ handles. This is the ideal setup to define a combinatorial version of knot Floer homology~\cite{OSknot,rasmussen}, known as \emph{grid homology}~\cite{MOS} (see also~\cite{SOS} for a detailed reference, and~\cite{BGH} for a lens space version). Despite its straightforward definition, a direct computation of grid homology becomes quickly unfeasible, as the chain complex has a number of generators which grows factorially in the size of the grid. It is however possible to greatly simplify the grid homology chain complex~\cite{beliakova}, and there is a relatively efficient algorithm to compute it \cite{baldwingillam} implemented in the software GridLink~\cite{gridlink}. Note that faster computation of knot Floer homology~\cite{OSknotfloerprogram} can be achieved using more recent techniques developed in~\cite{bordered_algebras}.

Thanks to their combinatorial definition, grid diagrams provide a convenient model to investigate on asymptotic properties of (bio)physical knots \cite{Witte, doig}. As an example, they have been used to quantify the intensity of crossing change-mediated fluxes between different knot types, and to demonstrate the dependency of those fluxes on the local geometry of knotted configurations. In turn, this provided strong evidence that the simplification action of specific DNA enzymes is driven by a geometric selection of sites \cite{griglie}. One further potential use of grid diagrams is to help with the search of band changes and the determination of Gordian distance between knot types~\cite{Moore,stolz, burnier, Cheston}. 

The most used program currently available to manipulate grid diagrams is GridLink. One of the main features of GridLink is a very efficient simplification function. Here we present GridPyM, a Sage~\cite{sage} compatible Python module that aims at enhancing certain functionalities present in GridLink, and at extending its manipulation capabilities. 

In particular, the module is built for generating large populations of complex grids to be employed for the statistical analysis of properties of the grid model. For this purpose, on one hand it is necessary to be able to quickly generate a large amount of random grid diagrams, possibly with specific properties (\emph{e.g.}~a fixed or bounded number of components). On the other hand, easy ways of efficiently simplifying large samples of grids are needed to accelerate invariants computations. 
Thus, GridPyM's focus is on generating, simplifying and randomising grids, rather than the computation of invariants, which can be carried out afterwards with efficient existing programs (such as~\cite{OSknotfloerprogram}, \cite{sage}, \cite{topoly}).

\section{Grid diagrams}

A \emph{grid diagram} $G$ is a $n \times n$ square plane grid, together with two sets of $n$ markings, conventionally denoted by $\XX$ and $\OO$.  Here $n \ge 2$ is called the \emph{grid number} or dimension of $G$. 

\begin{figure}[ht]
\centering
\captionsetup{width=13cm}
\includegraphics[width = 13cm]{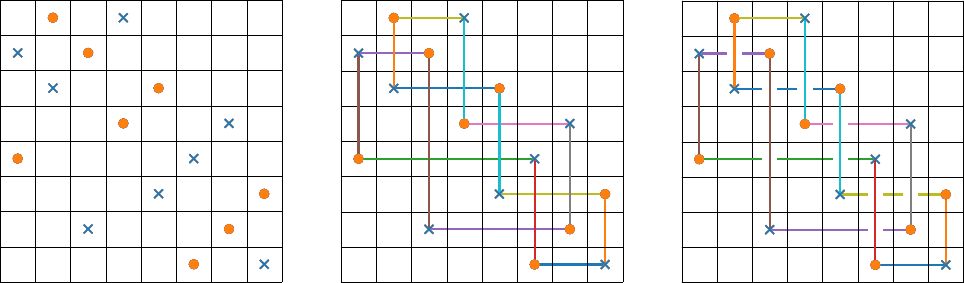}
\caption{From left to right: an example of a grid diagram, and how to obtain an oriented link diagram from it by resolving each double point as a vertical overpass.}
\label{fig:fromgrid2diagram}
\end{figure}

The grid is subdivided into $n^2$ little squares, and the markings are placed in these squares according to a ``sudoku'' rule, so that each row and column contains exactly one marking of each kind, and each square contains at most one marking. 
We adopt the convention according to which the markings are listed as ordered tuples $\XX = [X_1, \ldots, X_n]$ and $\OO = [O_1, \ldots, O_n]$, where the $i$-th component denotes the position (from the left, and starting from $0$) of the marking on the $i$-th row, where rows are enumerated from the bottom up. For example, the grid in the left of Figure~\ref{fig:fromgrid2diagram} is described by the pair of markings $\XX = [7,2,4,5,6,1,0,3]$ and $\OO = [5,6,7,0,3,4,1,2]$. Note that the combinatorial requirement on the markings implies that $\XX$ and $\OO$ are a pair of collision-free permutations of $\{0, \ldots, n-1\}$, where $n$ is the grid number.

To a grid diagram $G$ we can associate an oriented link diagram as follows: on each row of the grid, connect with a horizontal segment the unique $\OO$ marking to the only $\XX$, and do the same on columns from $\XX$ to $\OO$. 
Resolving each double point, with the convention that vertical lines are always overpasses, yields a diagram representing a link $L(G)$. A schematic outline of the process is shown in Figure~\ref{fig:fromgrid2diagram}. 

\begin{figure}[ht]
\centering
\captionsetup{width=13cm}
\includegraphics[width = 13cm]{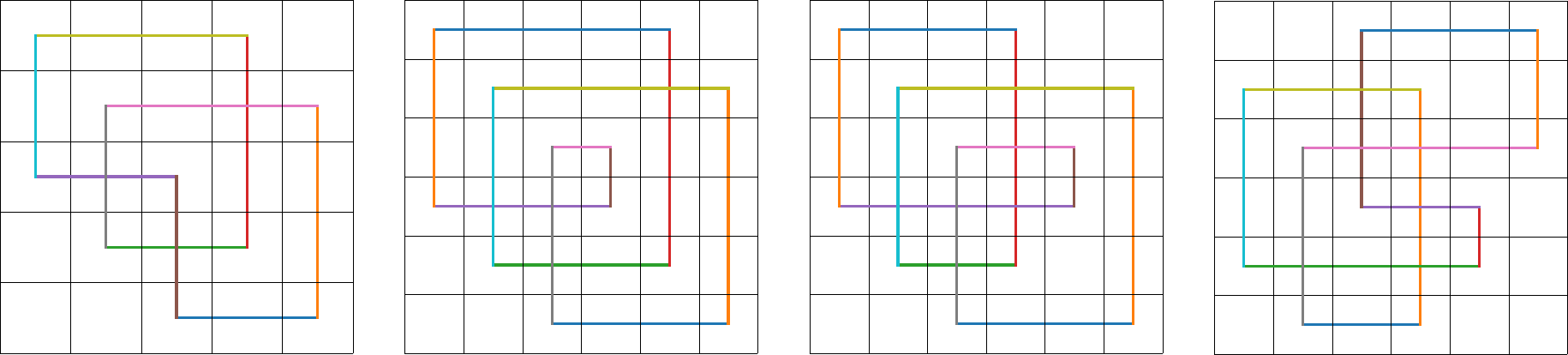}
\caption{From left to right: a grid diagram representing the (right) trefoil, the effect of applying a stabilisation, the effect of an interleaving commutation on the stabilised diagram and a shift along the vertical axis.}
\label{fig:cromwellmoves}
\end{figure}

It is easy to realise that any oriented link can be represented by infinitely many grid diagrams. Further, there is a finite set of moves, known as Cromwell moves (see Figure \ref{fig:cromwellmoves} for some examples), with the property that two grid diagrams are related by a finite sequence of such moves if and only if the links they represent are isotopic. These moves are therefore a grid-diagrammatic analogue of the classical Reidemeister moves on link diagrams.\\

\section{GridPyM: module description and main functions}
In what follows we give a brief description of the Python module GridPyM, and of some of the main functions implemented.

\subsection{Installation and use}

GridPyM is stored in the ``GridPythonModule'' GitHub repository (see Section~\ref{sec:availability}), together with some Jupyter notebooks~\cite{jupyter}, interactively explaining basic usage and functions. We tried to keep the module as self-contained as possible. The only external imports required are the rather standard libraries \url{sympy},  \url{random2} and \url{matplotlib}. As a consequence, there is full compatibility with Sage, and GridPyM can be imported and used in conjunction with Sage's built-in link functions. The module works for all versions of Python greater or equal to 3. The code will be periodically updated, and we refer to the GitHub repository for the most recent version. GridPyM can be downloaded from its GitHub repository. \\

To install GridPyM, write in a terminal 
\begin{verbatim}
    pip install GridPythonModule
\end{verbatim} 
To import GridPyM and its functions in a Python environment\slash terminal, write 

\begin{verbatim}
    import GridPythonModule
    from GridPythonModule import *
\end{verbatim}



\subsection{Grid generation}

GridPyM provides several different ways of generating grids, which are divided into:
\begin{itemize}
\item generating random grids with prescribed features (number of components and grid number) using the \texttt{generate\_random\_grid} commands,
\item generating a grid representing specific link types (\emph{e.g.}~torus links), 
\item loading from a library of low crossing number knots using \texttt{load\_knot},
\item loading grids representing low crossing number Legendrian knots (taken from the Legendrian knot atlas~\cite{ngatlas}), using \texttt{load\_legendrian\_knot}.
\end{itemize}

A list of available knot types can be accessed via the commands \texttt{available\_knots} and \texttt{available\_legendrian\_knots} respectively.

\begin{figure}
\centering
\captionsetup{width=12cm}
\includegraphics[width = 12cm]{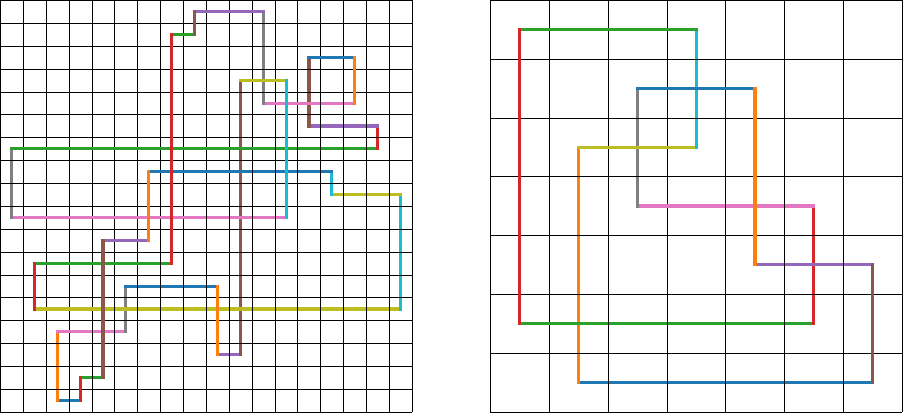}
\caption{The effect of the \texttt{destabilize\_all} function on a scrambled grid representing the $5_2$ knot. This function performs all possible ``generalised'' destabilisations on the grid.}
\label{fig:alldestab}
\end{figure}

\subsection{Grid handling and moves}

GridPyM includes all standard Cromwell moves (cyclic shifts, non-interleaving column\slash row commutations and (de)stabilisations, see Figure \ref{fig:cromwellmoves}). Rows and columns commutations might change the underlying link type of the diagram, depending on the local combinatorics. We include those resulting in a crossing changes~\cite[Prop~3.1.13]{SOS}, and band attachments. Band attachments are divided into \textit{coherent} or \textit{uncoherent}, depending on whether they change or not the number of components of the link. Examples on how to perform these moves are shown in the repository notebooks (see Section~\ref{sec:availability}).

\begin{figure}[ht]
\centering
\captionsetup{width=12cm}
\includegraphics[width = 12cm]{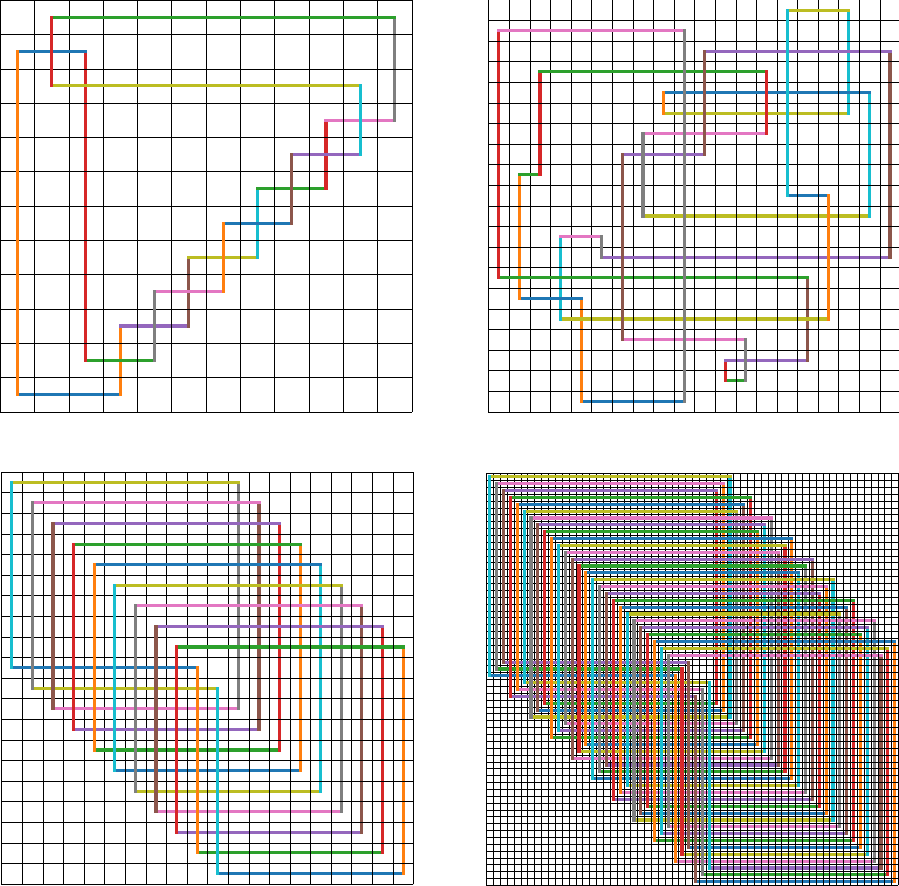}
\caption{Some examples of grid diagrams obtained with GridPyM's  drawing function \texttt{draw\_grid} (from top to bottom): the $8$ crossings negatively clasped twist knot, a random knot in grid number $20$, the $(11,9)$ torus knot and its flat $(3,1)$ cable. }
\label{fig:grids_from_code}
\end{figure}

It is also possible to perform some other operations to generate new links from given ones, such as disjoint union, connected sum, and some cables \emph{e.g.}~taking parallel copies of a knot (see Figure~\ref{fig:grids_from_code} for an example). Other standard knot-theoretic functions include inverting the orientation, taking the mirror image and rotating  the grid (note that rotating the grid by $\frac{\pi}{2}$ produces a representative for the mirror of the link). See the repository notebooks (Section~\ref{sec:availability}) for further examples.

\subsection{Grid simplification and randomisation}
The simplest way to reduce the size of a grid diagram is to perform a destabilisation. We implemented the function \texttt{destabilise\_all} that simplifies a grid by recursively destabilising all non-trivial configurations in which two markings are adjacent, as described in Figure~\ref{fig:alldestab}.

The function \texttt{destabilise\_all} is the key component of the main simplification function of GridPyM: \texttt{simplify\_grid}. The \texttt{simplify\_grid} function takes a grid diagram as input and tries (with customisable levels of effort) to simplify it, by performing destabilizations and a (customisable) number of random non-interleaving commutations and cyclic shifts. Note that by the \emph{monotonic simplification theorem}~\cite{dynnikov}, this process is (with probability $1$, in the absence of computational issues) guaranteed to return the $2\times2$ grid diagram, whenever $L(G)$ represents the unknot. See the repository notebooks (Section~\ref{sec:availability}) for examples. Further options include restricting to Cromwell moves that preserve the Legendrian or transverse class (so that the input and simplified grid represent the same Legendrian or transverse knot type). 

The function \texttt{scramble\_grid} provides the opposite operation: given a grid, it performs a customisable amount of random moves to it, making it on average more complex. Again, there is the option to preserve the Legendrian or transverse isotopy class.

\subsection{Grid and contact knot invariants}
As mentioned before, the suite of link invariants was kept to a minimum. We have thus only included invariants that can be efficiently computed from the grid, and that are not easily computable by other compatible programs. In particular, the only available topological invariant of the link type is the number of components. All other functions mentioned below are instead either invariants of the Legendrian or transverse class, or combinatorial invariants of the grid. In the former category we have the Thurston-Bennequin number, self-linking and rotation number. In the latter we have the grid number (so just the dimension of the given grid), crossing number, and the grid length; this is just the sum of the lengths of the vertical and horizontal segments forming the link.

\begin{figure}
\centering
\captionsetup{width=10cm}
\includegraphics[width=10cm]{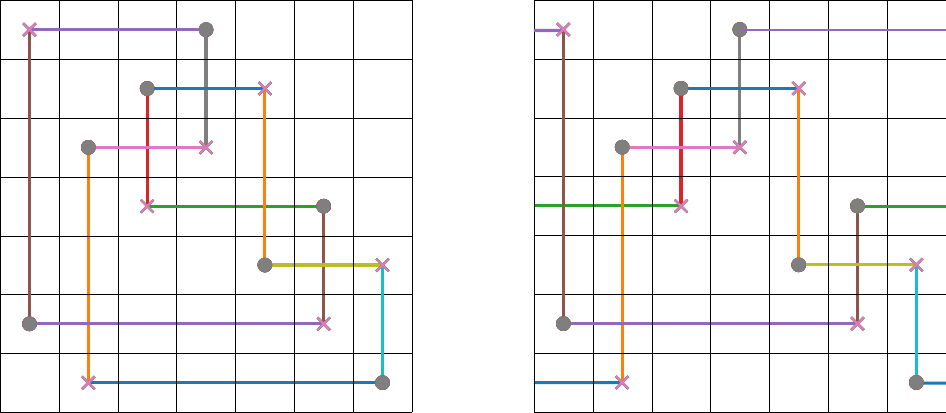}
\caption{A schematic representation showing one way to associate a braid diagram to a grid diagram.}
\label{fig:grid2braid}
\end{figure}

It is also possible to convert a given grid diagram into a braid (see Figure~\ref{fig:grid2braid}); this is done using the \url{convert_to_braid} function, whose output is a string containing the standard generators of a braid whose closure represents the link. For example, $[1,1,1]$ is the (right) trefoil, while $[1,-2,1,-2]$ is the figure-eight knot. Note that there is an option to automatically choose a pre-simplification function; in other words, the function first attempts to simplify the grid before computing the associated braid word.

Another way of passing from grid diagrams to knots is through Gauss codes. The output of the Gauss code and braid conversion functions are chosen to be compatible with Sage, as shown in the example below. (Note that the following code needs to be run on Sage.)

\begin{verbatim}
    >> G = load_knot(`7_3')
    >> braid_word = convert_to_braid(G)
    >> print(braid_word)
     [5, 4, 3, 2, 5, 4, 3, 2, 1, -2, 1]
    >> B = BraidGroup(6)
    >> K = Knot(B(braid_word))
    >> print(K.alexander_polynomial())
    2*t^-2 - 3*t^-1 + 3 - 3*t + 2*t^2
    >> g_code = Gauss_code(G)
    >> print(g_code)
    [[[1, 2, -3, 4, -2, -5, 6, -7, 8, -1, -4, 3, 5, -6, 7, -8]], 
    [1, -1, 1, 1, 1, 1, 1, 1]]
    >> L = Knot(g_code)
    >> show(L.plot())
    >> print(L.jones_polynomial())
    -t^9 + t^8 - 2*t^7 + 3*t^6 - 2*t^5 + 2*t^4 - t^3 + t^2
\end{verbatim}

\section{Sample computations}

In this final section, we collect some basic results obtained with GridPyM, mostly as a check of some of its functionalities. We analyse the distribution of several invariants for randomly generated grids in the grid number range $5-100$. For related results, see~\cite{griglie,doig}.

\begin{figure}[h!]
\centering
\captionsetup{width=13cm}
\includegraphics[width=13cm]{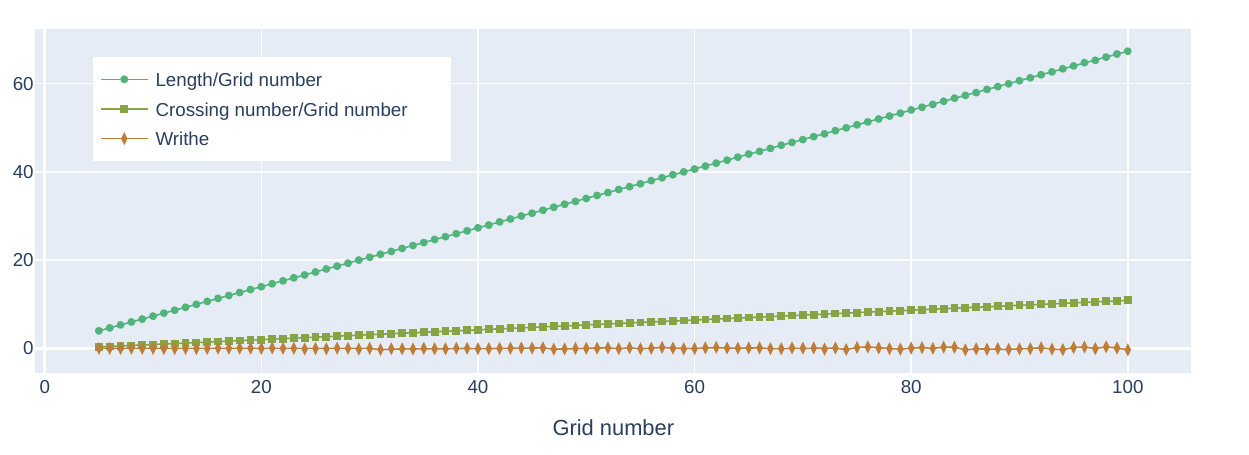}
\caption{A plot showing the growth of the average length of a grid divided by the grid's dimension, of the average crossing number divided by the grid's dimension and of the average writhe. As expected, the latter is distributed along the $y = 0$ line. On the other hand, both the crossing number and the length exhibit quadratic growth as functions of the grid number. The slope of the linear fit for length\slash grid number is $\sim0.6667$, and  $\sim0.1111$ for crossing number\slash grid number.}
\label{fig:cose}
\end{figure}
\begin{figure}
\centering
\captionsetup{width=13cm}
\includegraphics[width=13cm]{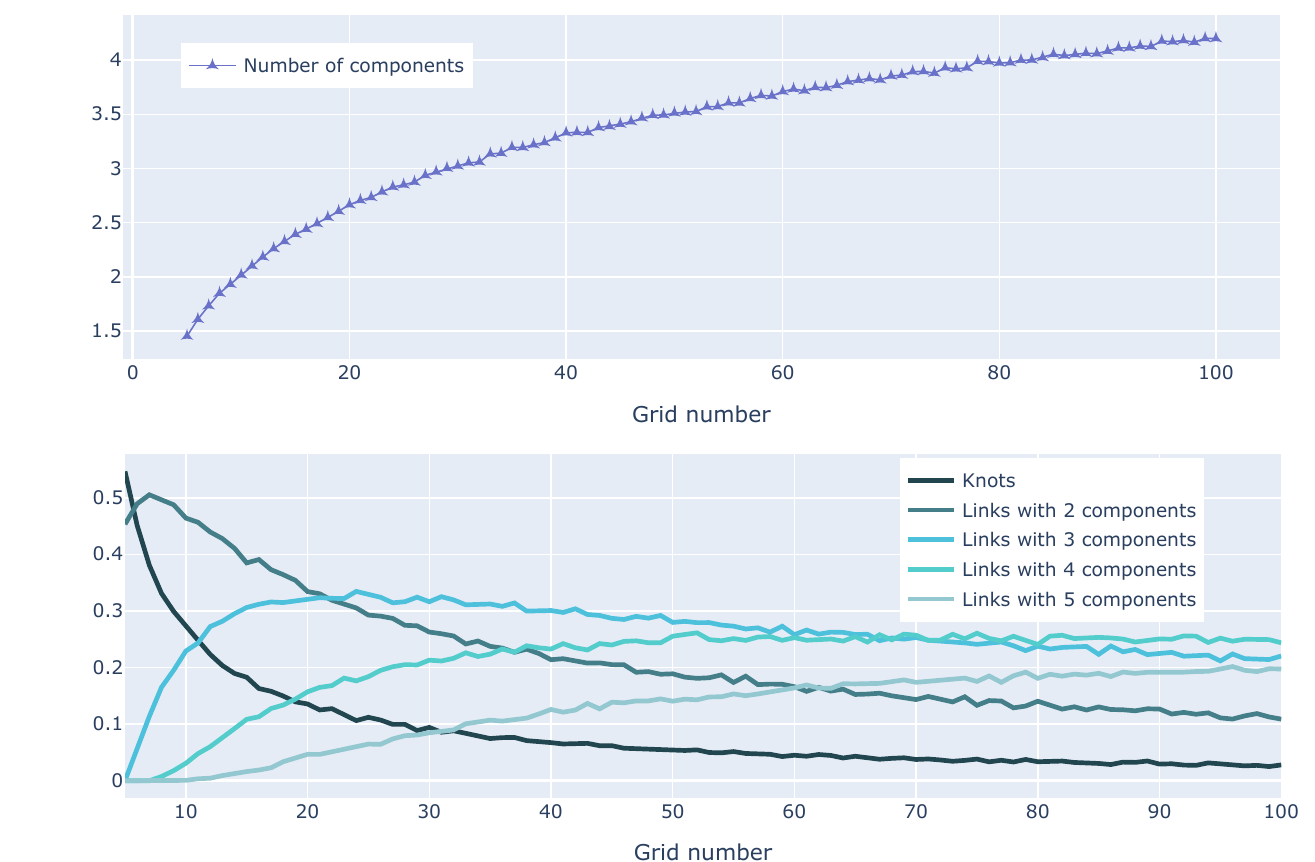}
\caption{In the top part of the figure, a plot of the average number of components in a sample population of random grid in the grid number range $5-100$. The bottom part of the figure shows the occurrence probabilities of grids with a given number of components as the grid number increases, for links with up to $5$ components. }
\label{fig:numb_comp}
\end{figure}
\begin{figure}
\centering
\captionsetup{width=13cm}
\includegraphics[width=13cm]{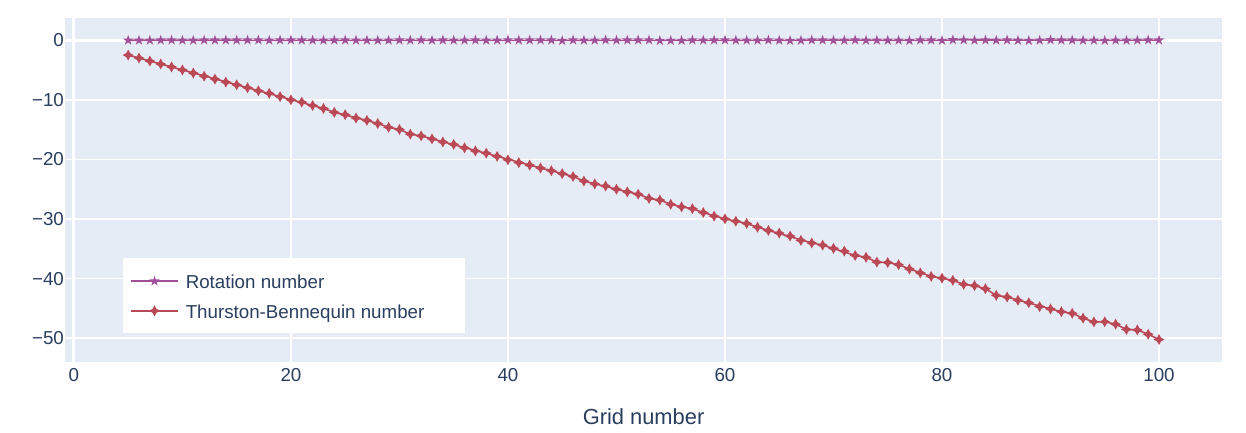}
\caption{The plot of grid contact invariants; consistently with their definition (\emph{cf}.~\cite{geiges}), the average rotation number is constantly $0$, while  the Thurston-Bennequin number is on a line with slope $-\frac{1}{2}$. }
\label{fig:average_contact}
\end{figure}

\section{Availability and implementation}\label{sec:availability}
GridPyM is freely available under the GNU General Public Licence v3.0. Source code and documentation can be found at \\
\begin{center}
    \url{https://github.com/agnesedaniele/GridPythonModule}
\end{center}

\subsection*{Acknowledgements} AB gratefully acknowledges funding through a MACSYS Centre Development initiative from the School of Mathematics \& Statistics, the Faculty of Science and the Deputy Vice Chancellor Research, University of Melbourne. DC was supported by Hodgson-Rubinstein’s ARC grant DP190102363 ``Classical And Quantum Invariants Of Low-Dimensional Manifolds''.

\bibliographystyle{amsplain}

\end{document}